\documentclass[12pt,leqno]{article}
\usepackage{amssymb,eufrak}
\begin{document}
\title{On L-Functions of Cyclotomic Function Fields}

\author{
Bruno Angl\`  es \\
Universit\'e de Caen, \\
Laboratoire Nicolas Oresme, CNRS UMR 6139, \\
Campus II, Boulevard Mar\'echal Juin, \\
BP 5186, 14032 Caen Cedex, France.\\
E-mail: angles@math.unicaen.fr}
\date{}

\maketitle
\begin{abstract} We study two criterions of cyclicity for divisor class groups of function fields, the first one involves Artin L-functions
and the second one involves ''affine'' class groups. We show that, in general, these two criterions are not
linked.
\end{abstract}\par
 Let $P$ be a prime of $\mathbb F_q[T]$ of degree $d$ and let $K_P$ be the $P$th cyclotomic function field. In this paper
we study the relation between the $p$-part of $Cl^0(K_P)$ and the zeta function of $K_P,$ where $p$ is the characteristic of $\mathbb
F_q.$\par Let $\chi$ be an even character of the Galois group of $K_P/\mathbb F_q(T),$ $\chi \not =1.$ Let $g(X,\overline{\chi })$ be the 
''congruent to one modulo $p$'' part of the L-function of $K_P/\mathbb F_q(T)$ associated to the character $\overline{\chi}.$ We have two
criterions of cyclicity (\cite{GOS}, chapter 8): if $\rm{deg}_Xg(X,\overline{\chi })\leq 1$ then $Cl^0(K_P)_p(\chi )$ is a cyclic $\mathbb
Z_{p}[\mu_{q^d-1}]$-module, and if $Cl(O_{K_P})_p(\chi)=\{ 0\}$ then  $Cl^0(K_P)_p(\chi )$ is a cyclic $\mathbb
Z_{p}[\mu_{q^d-1}]$-module. David Goss has obtained   that if $Cl(O_{K_P})_p(\chi)$ is trivial  
then 
$g(X,\overline{\chi })$ is of degree at most one  (\cite{GOS}, Theorem 8.21.2).  Unfortunately,
there is a gap in the proof of this result. In fact,  we show that in general  
$Cl(O_{K_P})_p(\chi)=\{ 0\}$ does not imply  $\rm{deg}_Xg(X,\overline{\chi })\leq 1$ (Proposition
\ref{Proposition3}). We also prove that if $i$ is a $q$-magic number and if $\omega_P$ is the
Teichm\"uller character at
$P,$ then
$g(X,\omega_P^i)$ has simple roots when $i\equiv 0\pmod{q-1}$ (Proposition
\ref{Proposition7}).\par Note that  Goss conjectures that if $i$ is a $q$-magic number then
$\rm{deg}_Xg(X,\omega_P^i)\leq 1.$ This problem is still open and can be viewed as an analogue of
Vandiver's Conjecture for function fields (see section 5).\par
${}$\par
The author thanks David Goss (the proof of Lemma \ref{Lemma7} was communicated to the author by David Goss) and Philippe Satg\'e for several
fruitfull discussions.\par

\section{Notations}\par
${}$\par
Let $\mathbb F_q$ be a finite field having $q$ elements, $q=p^s$ where
$p$ is the characteristic of $\mathbb F_q.$ Let $T$ be an indeterminate
over $\mathbb F_q$ and set $A=\mathbb F_q[T],$ $k=\mathbb F_q(T).$ We
denote the set of monic elements of $A$ by $A^+.$ A prime of $A$ is a
monic irreducible polynomial in $A.$ We fix $\overline{k}$ an algebraic
closure of $k.$ We denote the unique place of $k$ which is a pole of $T$
by $\infty .$\par
${}$\par
Let $L/k$ be a finite geometric  extension of $k,$ $L\subset \overline{k}.$ We
set:\par
\noindent - $O_L:$ the integral closure of $A$ in $L,$\par
\noindent - $O_L^*:$ the group of units of $O_L,$\par
\noindent - $S_{\infty}(L):$ the set of places of $L$ above $\infty ,$\par
\noindent - $Cl^0(L):$ the group of divisors of degree zero of $L$ modulo
the group of principal divisors,\par
\noindent - $Cl(O_L):$ the ideal class group of $O_L,$\par
\noindent - $R(L):$ the groupe of divisors of degree zero with supports in $S_{\infty}(L)$ modulo the group of principal divisors
with supports in $S_{\infty}(L).$ \par
If $d$ is the greatest common divisor of the degrees of the elements in $S_{\infty}(L),$ we have the following exact sequence:
$$0\rightarrow R(L)\rightarrow  Cl^0(L)\rightarrow Cl (O_L) \rightarrow \frac{\mathbb Z}{d\mathbb Z}\rightarrow 0\, .$$
\par

Let $P$ be a prime of $A$ of degree $d.$ We denote the $P$th cyclotomic
function field by $K_P$ (see \cite{GOS}, chapter 7, and \cite{HAY}).
Recall that $K_P/k$ is the maximal abelian extension of $k$ contained in
$\overline{k}$ such that:\par
\noindent - $K_P/k$ is unramified outside of $P, \infty ,$\par
\noindent - $K_P/k$ is tamely ramified at $P,\infty ,$\par
\noindent - for every place $v$ of $K_P$ above $\infty ,$ the completion
of $K_P$ at $v$ is equal to $\mathbb
F_q((\frac{1}{T}))(^{q-1}\sqrt{-T}).$\par
We recall that $\rm{Gal}(K_P/k)\simeq (A/PA)^*,$ and that the
decomposition group of $\infty$ in $K_P/k$ is equal to its inertia group
and is isomorphic to $\mathbb F_q^*.$\par
${}$\par
Let $E/\mathbb F_q$ be a global function field and let $F/E$ be a finite
geometric abelian extension. Set $G=\rm{Gal}(F/E)$ and $\widehat
{G}=\rm{Hom}(G,\mathbb C^*).$\par
 Let $\chi \in \widehat {G},$ $\chi \not
=1,$ we set:
$$L(X,\chi )=\prod_{v\, \rm{place\, of}\, E}(1-\chi
(v)X^{\rm{deg}v})^{-1},$$
Where $\chi (v)=0$ if $v$ is ramified in $F^{\rm{Ker}(\chi)}/E,$ and if
$v$ is unramified in $F^{\rm{Ker}(\chi)}/E,$ $\chi(v)=\chi
((v,F^{\rm{Ker}(\chi)}/E)),$ where  $(.,F^{\rm{Ker}(\chi)}/E)$ is the
global reciprocity map. If $\chi =1,$ we set $L(X,\chi )=L_E(X)$ where
$L_E(X)$ is the numerator of the zeta function of $E.$\par
Therefore, if $L_F(X)$ is the numerator of the zeta function of $F,$ we
get:
$$L_F(X)=\prod_{\chi \in \widehat{G}}L(X,\chi )\, .$$\par
Let $\Delta $ be a finite abelian group and let $M$ be a $\Delta $-module. Let $\ell$ be a primÿe
number such that $\mid \Delta \mid \not \equiv 0\pmod{\ell }.$ We fix an embedding of
$\overline{\mathbb Q}$ in $\overline{\mathbb Q_{\ell}}.$  Let
$W=\mathbb Z_{\ell }[\mu_{\mid \Delta \mid}].$ For $\chi \in \widehat{\Delta},$ we set:
$$e_{\chi}=\frac{1}{\mid \Delta \mid}\sum_{\delta \in \Delta}\chi(\delta)\delta^{-1}\in
W[\Delta],$$
and:
$$M_{\ell}(\chi )=e_{\chi}(M\otimes_{\mathbb Z}W).$$
Thus, we have:
$$ M\otimes_{\mathbb Z}W =\bigoplus_{\chi \in \widehat{\Delta}}M_{\ell}(\chi )\, .$$

\section{Cyclotomic Function Fields and Artin-Schreier Extensions}
${}$\par
Let $Q$ be a prime of $A$ of degree $n,$ write $Q(T)=T^n+\alpha T^{n-1}+\cdots ,$ $\alpha \in \mathbb F_q.$ We
set: $i(Q)=Tr_{\mathbb F_q/\mathbb F_p}(\alpha ).$ Let $a\in A,$ $a\not =0,$ we set:
$$i(a)=\sum_{Q\, \rm{prime \, of}\, A}v_Q(a)i(Q)\in \mathbb F_p,$$
where $v_Q$ is the normalized $Q$-adic valuation on $k.$\par
${}$\par
Let $\theta \in\overline{k}$ such that $\theta^p-\theta =T.$ Set $\widetilde {A}=\mathbb F_q[\theta ],$
$\widetilde{k}=\mathbb F_q(\theta)$ and $G=\rm{Gal}(\widetilde {k}/k).$ Note that $\widetilde{k}/k$ is
unramified outside $\infty $ and totally ramified at $\infty .$ Let $\widetilde{\infty }$ be the unique place
of $\widetilde {k}$ above $\infty .$\par
\newtheorem{Lemma1}{Lemma}[section]
\begin{Lemma1} \label{Lemma1}
Let $(.,\widetilde {k}/k)$ be the usual Artin symbol. For $a\in A\setminus \{ 0 \} :$
$$(a,\widetilde{k}/k)(\theta )=\theta -i(a).$$
\end{Lemma1}
\noindent{\sl Proof} By the classical properties of the Artin symbol, it is enough to prove the Lemma when $a$
is a prime of $A.$ Thus, let $P$ be a prime of $A$ of degree $d.$ We have:
$$(P,\widetilde{k}/k)(\theta )\equiv \theta ^{q^d}\pmod{P}.$$
But, for $n\geq 0,$ we have:
$$\theta^{p^n}=\theta +T+T^p+\cdots +T^{p^{n-1}}.$$
Therefore:
$$\theta^{q^d}\equiv \theta -i(P)\pmod{P}.$$
The Lemma follows. $\diamondsuit$\par
\newtheorem{Lemma2}[Lemma1]{Lemma}
\begin{Lemma2} \label{Lemma2}
Let $P$ be a prime of $A$ of degree $d$ such that $i(P)\not =0.$ Then $P$ is a prime of $\widetilde {A}$ of
degree $pd.$ Let $\widetilde{K_P}$ be the $P$th cyclotomic function field for the ring $\widetilde{A},$ then
$K_P\subset \widetilde{K_P}.$
\end{Lemma2}
\noindent{\sl Proof}  We have $-T=-\theta^p(1-\theta^{1-p}).$ Note that:
$$1-\theta^{1-p}\in \mathbb (F_q((\frac{1}{\theta }))^*)^{q-1}.$$
Therefore:
$$^{q-1}\sqrt{-T}\in F_q((\frac{1}{\theta }))(^{q-1}\sqrt{-\theta }).$$
Thus:\par
\noindent - $\widetilde{k}K_P/\widetilde {k}$ is unramified outside $P, \widetilde{\infty },$\par
\noindent - $\widetilde{k}K_P/\widetilde {k}$ is tamely ramified at $P, \widetilde {\infty},$ \par
\noindent - for every place $w$ of $\widetilde{k}K_P$ above $\widetilde{\infty },$ the completion of 
$\widetilde{k}K_P$  at $w$ is contained in  $F_q((\frac{1}{\theta }))(^{q-1}\sqrt{-\theta }).$\par
The Lemma follows by class field theory. $\diamondsuit$\par
${}$\par
Let $P$ be a prime of $A,$ $\rm{deg}_TP(T)=d$ and $i(P)\not =0.$ Let $L=\widetilde{k}K_P\subset
\widetilde{K_P}.$ Let $\Delta=\rm{Gal}(K_P/k)\simeq \rm{Gal}(L/\widetilde{k}).$ We have an isomorphism
compatible to class field theory: $\widehat {\Delta}\rightarrow \widehat {\rm{Gal}(L/\widetilde {k})},$ $\chi
\mapsto \widetilde {\chi}=\chi \circ N_{\widetilde{k}/k}.$ We fix $\zeta_p\in \overline{\mathbb Q}$ a primitive
$p$th root of unity.
\newtheorem{Lemma3}[Lemma1]{Lemma}
\begin{Lemma3}\label{Lemma3}
${}$ \par
\noindent (1) Let $\chi \in \widehat{\Delta},$ $\chi \not =1.$ Let $L(X,\widetilde{\chi})$ be the Artin L-function
relative to $L/\widetilde{k}$ and to the character $\widetilde{\chi }.$ We have:
$$L(X,\widetilde{\chi })=\prod_{\phi \in \widehat{G}}L(X,\phi \chi ),$$
where $L(X,\phi \chi)$ is the Artin L-function relative to $L/k$ and the character $\phi \chi .$\par
\noindent (2) Let $\chi \in \widehat{\Delta},$ $\chi \not =1,$ $\chi$ even (i.e. $\chi (\mathbb F_q^*)=\{ 1\}$). Then:
$$\frac{L(X,\widetilde{\chi})}{L(X,\chi)}\equiv (1-X)^{p-1}L(X,\chi)^{p-1}\pmod{(1-\zeta_p)}.$$
\end{Lemma3}
\noindent{\sl Proof} Te assertion (1) is a consequence of the usual properties of Artin L-functions. Now, let $\phi \in
\widehat{G},$ $\phi \not =1.$ Since $\phi \chi$ is ramified at $\infty,$ we get:
$$L(X,\phi \chi )=\sum_{n\geq 0}(\sum_{a\in A^+,\, \rm{deg}(a)=n}\phi (a)\chi(a))X^n.$$
Thus:
$$L(X,\phi \chi)\equiv \sum_{n\geq 0}(\sum_{a\in A^+,\, \rm{deg}(a)=n}\chi (a)))X^n\pmod{(1-\zeta_p)}.$$
But, since $\chi$ is even, we have $\chi (\infty )=1.$  Therefore:
$$L(X,\phi \chi)\equiv (1-X)L(X,\chi)\pmod{(1-\zeta_p)}.$$
The Lemma follows. $\diamondsuit$\par
let $i\in\mathbb F_p$ and let $\sigma_i\in G$ such that $\sigma_i(\theta )=\theta -i.$ Let $\psi \in \widehat {G}$ given
by $\psi(\sigma_i)=\zeta_p^i.$
\newtheorem{Lemma4}[Lemma1]{Lemma}
\begin{Lemma4}\label{Lemma4}
Let $\chi \in \widehat{\Delta},$ $\chi$ even and non-trivial.\par
\noindent (1) Let $\phi \in \widehat{G},$ $\phi \not =1.$ Let $\sigma \in \rm{Gal}(\mathbb Q(\zeta_p)/\mathbb
Q)$ such that $\phi =\psi^{\sigma}.$ Then:
$$L(X,\phi \chi)=L(X,\psi \chi )^{\sigma}.$$
Furthermore $\rm{deg}_XL(X,\phi \chi )=d.$\par
\noindent (2) We have:
$$L(1,\psi \chi)\equiv (\sum_{a\in A^+, \rm{deg}(a)\leq d}i(a)\chi (a))(\zeta_p-1)\pmod{(1-\zeta_p)^2}.$$
\end{Lemma4}
\noindent{\sl Proof} Let $\mathbb Q (\chi )$ be the abelian extension of $\mathbb Q$ obtained by adjoining to
$\mathbb Q$ the values of $\chi .$ Let $\mathbb Z[\chi ]$ be the ring of integers of $\mathbb Q( \chi ).$ Note 
that $p$ is unramified in $\mathbb Q(\chi )$ and:
$$\rm{Gal}(\mathbb Q (\chi )(\zeta_p)/\mathbb Q(\chi ))\simeq \rm{Gal}(\mathbb Q(\zeta_p)/\mathbb Q).$$
Since $L(X,\phi \chi)$ is a polynomial in $\mathbb Z[\chi ][\zeta_p][X],$ we have:
$$L(X,\phi \chi)=L(X,\psi \chi )^{\sigma}.$$
Since $\chi$ and $\widetilde \chi$ are non-trivial even characters, we have:
$$\rm{deg}_XL(X,\widetilde {\chi})=pd-2,$$
and:
$$\rm{deg}_XL(X,\chi)=d-2.$$
Therefore $\rm{deg}_XL(X,\phi \chi)=d.$\par
Now, we have:
$$L(X,\psi \chi)=\sum_{n=0}^d(\sum_{a\in A^+\, \rm{deg}(a)=n}\zeta_p^{i(a)}\chi (a) )X^n.$$
But recall that:
$$\zeta_p^{i(a)}\equiv 1+i(a)(\zeta_p-1)\pmod{(1-\zeta_p)^2}.$$
Thus, since $\chi$ is even and  non-trivial, we get:
$$L(X,\psi \chi) \equiv L(X,\chi )(1-X) +(\zeta_p -1)(\sum_{n=1}^d(\sum_{a\in A^+\, \rm{deg}(a)=n}i(a)\chi (a))
X^n) \pmod{(1-\zeta_p)^2}.$$
The Lemma follows. $\diamondsuit$\par
We are now ready to prove the main result of this section:
\newtheorem{Proposition1}[Lemma1]{Proposition}
\begin{Proposition1}  \label{Proposition1}
Let $\chi \in \widehat {\Delta},$ $\chi \not =1,$ $\chi$ even. Let $W=\mathbb Z_p[\mu_{q^d-1}].$  We have:
$$\rm{Long}_W(\frac{Cl(O_L)_p(\widetilde {\chi})}{Cl(O_{K_P})_p(\chi )})\geq 1 \Leftrightarrow \sum_{a\in A^+\, \rm{deg}(a)\leq d}
i(a)\overline{\chi }(a) \equiv 0\pmod{p}.$$
\end{Proposition1}
\noindent{\sl Proof} Fix $\tau$ a generator of $G\simeq \rm{Gal}(L/K_P).$ Let $\varepsilon \in O_L^*.$ Since $L/K_P$ is totally
ramified at any prime above $\infty,$ there exists $\zeta \in \mathbb F_q^*$ such that $\tau (\varepsilon) =\zeta \varepsilon .$
But $\tau^p(\varepsilon)=\zeta^p\varepsilon =\varepsilon.$ Since we are in characteristic $p,$ we deduce that $\varepsilon \in
O_{K_P}^*.$ Therefore:
$$O_L^*=O_{K_P}^*.$$
Let $I$ be an ideal of $O_{K_P}$ such that $IO_L=\alpha O_L$ for some $\alpha \in O_L.$ Then, there exists $\varepsilon \in O_L^*$
 such that $\tau (\alpha )= \varepsilon \alpha .$ Since $O_L^*=O_{K_P}^*$ and since $\tau$ is of order $p,$ we deduce that $\alpha
\in O_{K_P}.$ This implies that:
$$Cl(O_{K_P})\hookrightarrow Cl(O_L).$$
One can also show that:
$$Cl^0(K_P)\hookrightarrow Cl^0(L).$$
Set $\Delta ^+=\frac{\Delta }{\mathbb F_q^*}.$  Let $\cal I$ be the augmentation ideal of $\mathbb F_p[\Delta^+].$ One sees that 
we have the following isomorphism of $\Delta $-modules:
$$\frac{R(L)}{R(K_P)}\otimes_{\mathbb Z}\mathbb Z_p\simeq \cal I.$$
This implie that we have the following exact sequence of $W$-modules:
$$0\rightarrow \frac{W}{pW}\rightarrow \frac{Cl^0(L)_p(\widetilde {\chi })}{Cl^0(K_P)_p(\chi )}\rightarrow 
\frac{Cl(O_L)_p(\widetilde{\chi })}{Cl(O_{K_P})_p(\chi )}\rightarrow 0\, .$$
Now, by the results of Goss and Sinnott (\cite{GOS&SIN}):
$$\rm{Long}_WCl^0(L)_p(\widetilde{\chi })= v_p(L(1,\overline{\widetilde{\chi }})),$$
and
$$\rm{Long}_WCl^0(K_P)_p(\chi )=v_p(L(1,\overline{\chi })).$$
Thus by Lemma \ref{Lemma3}:
$$\rm{Long}_W(\frac{Cl(O_L)_p(\widetilde {\chi})}{Cl(O_{K_P})_p(\chi )})=(p-1)v_p(L(1,\psi \overline{\chi }))-1.$$
It remains to apply Lemma \ref{Lemma4}. $\diamondsuit$\par
\section{Derivatives of L-functions}
${}$\par
Let $P$ be a prime of $A$ of degree $d.$ We fix an embedding of $\overline{\mathbb Q}$ in $\overline{\mathbb Q_p}.$ Set $\Delta
=\rm{Gal}(K_P/k)$ anf
$W=\mathbb Z_p[\mu_{q^d-1}].$ We fix an isomorphism
$\Phi_P:\, A/PA\rightarrow W/pW.$ Then $\Phi_P$ induces an isomorphism:
$$\omega_P:\, \Delta \rightarrow \mu_{q^d-1}\subset W^*.$$
The morphism $\omega_P$ is called  ''the'' Teichm\" uller character at $P.$
Note that $\widehat{\Delta }$ is a cyclic group and $\omega_P$ is a generator of this group.\par
Let $i\in \mathbb N,$ set:\par
\noindent - $\beta (0)=1,$\par
\noindent - $\beta (i)=\sum_{a\in A^+}a^i$ if $i\geq 1,$ $i\not \equiv 0\pmod{q-1},$\par
\noindent - $\beta (i)=-\sum_{a\in A^+}\rm{deg}(a)a^i$if $i\geq 1,$ $i\equiv 0\pmod{q-1}.$\par
\noindent One can prove that for all $i\in \mathbb N,$ $\beta (i) \in A.$ We also see that:
$$\forall i\in \mathbb N,\, 0\leq i\leq q^d-2,\, \Phi_P(\beta (i))\equiv L(1, \omega_P^i) \pmod{p}.$$
Therefore, if $1\leq i \leq q^d-2,$ by the results of Goss and Sinnott (\cite{GOS&SIN}), we have:
$$\rm{Long}_WCl^0(K_P)_p(\omega_P^{-i})\geq 1 \Leftrightarrow \beta(i) \equiv 0\pmod{P}.$$
The numbers $\beta(i)$ are called the Bernoulli-Goss polynomials.\par
 Recall that we have a surjective morphism of $\Delta$-modules:
$$W[\Delta^+]\rightarrow R(K_P)\otimes_{\mathbb Z}W,$$
where $\Delta^+=\Delta /\mathbb F_q^*.$ Thus for $\chi \in \widehat{\Delta},$ $\chi$ even, $R(K_P)_p(\chi)$ is a cyclic $W$-module.
But,  for such a character, we have the exact sequence of $W$-modules:
$$0\rightarrow R(K_P)_p(\chi )\rightarrow Cl^0(K_P)_p(\chi) \rightarrow Cl(O_{K_P})_p(\chi )\rightarrow 0\, .$$
This implies that, if $Cl(O_{K_P})_p(\chi )=\{ 0\},$ $Cl^0(K_P)_p(\chi )$ is a cyclic $W$-module.\par
David Goss has shown (\cite{GOS}, Corollary 8.16.2) that for  $\chi$ is even, $\chi \not =1,$  if $L'(1,\overline{\chi })\not
\equiv 0 \pmod{p}$ (here $L'(1,\overline{\chi})$ is the derivative of $L(X,\overline{\chi })$ taken at $X=1$), then
$Cl^0(K_P)_p(\chi )$ is a cyclic $W$-module.\par
Therefore a natural question arise. Let $\chi \in \widehat{\Delta},$ $\chi \not =1,$ $\chi$ even. Assume that $L(1,\overline{\chi
})\equiv 0 \pmod{p}.$ Do we have: 
$$Cl(O_{K_P})_p(\chi)=\{ 0\} \Rightarrow L'(1, \overline{\chi })\not \equiv 0 \pmod{p}?$$
Our aim in this section is to show that in general the answer is no.\par
${}$\par
Let $d$ be an integer, $d\geq 1.$ For $i\in \{ 1,\cdots, q^d-2\},$ we set:
$$\gamma (d,i)=\sum_{a\in A^+,\, \rm{deg}(a)\leq d}i(a)a^i.$$
\newtheorem{Lemma5}{Lemma}[section]
\begin{Lemma5} \label{Lemma5}
Let $\tau \in \rm{Gal}(\mathbb F_q(T)/\mathbb F_q(T^p-T))$ such that $\tau (T)=T+1.$ Let $i\in \{ 1,\cdots ,q^d-2\},$ $i\equiv 0\pmod{q-1}.$
Recall that $q=p^s.$ We have:
$$\tau(\gamma (d,i))=\gamma(d,i)+s\beta (i).$$
\end{Lemma5}
\noindent{\sl Proof} Let $Q$ be a prime of $A$ of degree $n.$ Write $Q=T^n+\alpha T^{n-1}+\cdots,$ where $\alpha \in \mathbb F_q.$ Then $\tau
(Q)=T^n+(\alpha +n)T^{n-1}+\cdots$ . Therefore $i(\tau (Q))=i(Q)+s\rm{deg}(Q).$ This implies that:
$$\forall a\in A\setminus \{ 0\},\, i(\tau (a))=i(a)+s\rm{deg}(a).$$
Now:
$$\tau (\gamma (d,i))=\sum_{a\in A^+,\, \rm{deg}(a)\leq d}i(a) \tau (a)^i.$$
Therefore:
$$\tau (\gamma (d,i))=\sum_{a\in A^+,\, \rm{deg}(a)\leq d}(i(\tau(a))-s\rm{deg}(a))\tau (a)^i.$$
Thus:
$$\tau(\gamma(d,i))=\sum_{a\in A^+,\, \rm{deg}(a)\leq d}i(\tau (a))\tau (a)^i -s\sum_{a\in A^+,\, \rm{deg}(a)\leq d}\rm{deg}(\tau (a)) \tau
(a)^i.$$
Observe that $\sum_{a\in A^+,\, \rm{deg}(a)\leq d}i(\tau (a))\tau (a)^i =\gamma (d,i)$ and $-\sum_{a\in A^+,\, \rm{deg}(a)\leq d}\rm{deg}(\tau
(a)) \tau (a)^i=\beta (i).$ $\diamondsuit$\par
\newtheorem{Proposition2}[Lemma5]{Proposition}
\begin{Proposition2} \label{Proposition2}
Let $P$ be a prime of $A$ of degree $d$ such that $i(P)\not =0.$ Set $Q(T)=P(T^p-T).$ Then $Q$ is a prime of $A$ of degree $pd.$ Let $i$ be
an integer such that $1\leq i\leq q^d-2,$ $i\equiv 0 \pmod{q-1}$ and $Cl(O_{K_P})_p(\omega_P^{-i})=\{ 0\} .$ Then:
$$\rm{Long}_W Cl(O_{K_Q})_p(\omega_Q^{-i(q^{pd}-1)/(q^d -1)})\geq  1 \Leftrightarrow \gamma (d,i)\equiv 0 \pmod{P}.$$
\end{Proposition2}
\noindent{\sl Proof} We have:
$$\Phi_P(\gamma (d,i))\equiv \sum_{a\in A^+,\, \rm{deg}(a)\leq d}i(a) \omega_P^i(a)\pmod{p}.$$
It remains to apply Proposition \ref{Proposition1}. $\diamondsuit$ \par
\newtheorem{Lemma6}[Lemma5]{Lemma}
\begin{Lemma6} \label{Lemma6}
Assume $p\not = 2.$ Let $d\geq 1$ be an integer.  There exists a prime $P$ in $A,$ $\rm{deg}(P)=d,$  such that
$i(P(T))i(P(T+1))\not =0.$
\end{Lemma6}
\noindent{\sl Proof} Let $Q$ be a prime of $A$ of degree $d$ such that $i(Q)\not = 0.$ Such a prime exists by the normal basis Theorem. Fix
$\overline{\mathbb F_q}$ an algebraix closure of $\mathbb F_q.$ We assume that $i(Q(T+1))=0.$ Write $Q=T^d+\alpha T^{d-1}+\cdots$ . Then 
$Tr_{\mathbb F_q/\mathbb F_p}(\alpha)=-sd.$ Therefore $sd\not \equiv 0 \pmod{p}.$Let $\theta \in\overline{\mathbb F_q}$ such that $Q(\theta
)=0.$ We observe that:
$$\forall \zeta \in \mathbb F_p,\, Tr_{\mathbb F_{q^d}/\mathbb F_p}(\zeta \theta )=-\zeta sd.$$
Since $p\geq 3,$ we can find $\zeta \in \mathbb F_p^*$ such that $-\zeta sd \not = -sd.$ Set $P(T)=\rm{Irr}(\zeta \theta , \mathbb F_q ; T).$
Then $P$ is a prime of degree $d$ such that $i(P)i(\tau (P))\not = 0.$ $\diamondsuit$ \par
\newtheorem{Proposition3}[Lemma5]{Proposition}
\begin{Proposition3} \label{Proposition3}
Assume that $p\not =2$ and $s\not \equiv 0\pmod{p}.$ Let $d$ be an integer, $d\geq 2,$ and let $P$ be a prime of degree $d$ such that
$i(P(T))i(P(T+1))\not =0.$ Set $Q(T)=P(T^p-T).$ Then:\par
\noindent - $L(1,\omega_Q^{-(q-1)(q^{pd}-1)/(q^d-1)})\equiv 0\pmod{p},$\par
\noindent - $L'(1,\omega_Q^{-(q-1)(q^{pd}-1)/(q^d-1)})\equiv 0\pmod{p},$\par
\noindent - $Cl(O_{K_Q})_p(\omega_Q^{-(q-1)(q^{pd}-1)/(q^d-1)})=\{ 0\}.$
\end{Proposition3}
\noindent{\sl Proof} Set $R=P(T+1)$ and $Z=R(T^p-T).$ We observe that we have an isomorphism:
$$Cl(O_{K_Q})_p(\omega_Q^{-(q-1)(q^{pd}-1)/(q^d-1)})\simeq Cl(O_{K_Z})_p(\omega_Z^{-(q-1)(q^{pd}-1)/(q^d-1)}).$$
Not also that $\beta (q-1)=1.$ Thus:
$$Cl(O_{K_P})_p(\omega_P^{-(q-1)})=Cl(O_{K_R})_p(\omega_R^{-(q-1)})=\{ 0\}.$$
We have:
$$L(1,\omega_Q^{-(q-1)(q^{pd}-1)/(q^d-1)})\equiv L(1,\omega_Z^{-(q-1)(q^{pd}-1)/(q^d-1)})\equiv 0\pmod{p}.$$
And, by Lemma \ref{Lemma3}, since $p\geq 3$: 
$$L'(1,\omega_Q^{-(q-1)(q^{pd}-1)/(q^d-1)})\equiv L'(1,\omega_Z^{-(q-1)(q^{pd}-1)/(q^d-1)})\equiv 0\pmod{p}.$$
Suppose that we have $Cl(O_{K_Q})_p(\omega_Q^{-(q-1)(q^{pd}-1)/(q^d-1)})\not =\{ 0\}.$ Then by Proposition\ref{Proposition2}:
$$\gamma(d,q-1)\equiv 0\pmod{P},$$
 and also:
$$\gamma (d,q-1)\equiv 0\pmod{R}.$$
Thus:
$$\tau (\gamma (d,q-1))\equiv 0\pmod {\tau (P)}.$$
Now, by Lemma \ref{Lemma5}, and the fact that $\tau (P)=R,$  we get:
$$\gamma(d,q-1)+s\beta (q-1)\equiv 0\pmod{R}.$$
Therefore we get $s\equiv 0\pmod{p}$ which is a contradiction. The Proposition follows. $\diamondsuit$ \par
\section{Cyclicity of Class Groups and L-Functions}
${}$\par
Let $E/\mathbb F_q$ be a global function field and let $F/E$ be a finite geometric abelian extension. Set $\Delta =\rm{Gal}(F/E).$ Let $\ell$
be a prime number. Let's recall some well-known facts about $L$-functions.\par
${}$\par
Set $T_{\ell }=\rm{Hom}(\mathbb Q_{\ell }/\mathbb Z_{\ell } , J)$ where  $J$ is the inductive limit of the $ Cl^0(\mathbb F_{q^n}F),$ $n\geq
1.$ We fix an embedding of $\overline{\mathbb Q}$ in $\overline{\mathbb Q_{\ell}}.$ Let $\gamma$ be the Frobenius of $\mathbb F_q.$ Then
$\gamma$ and $\Delta$ act on $T_{\ell }.$\par
If $\ell \not =p,$ we have (see \cite{ROS},chapter 15):
$$\rm{Det}(1-\gamma X\mid_{T_{\ell}})=L_F (X),$$
 where $L_F (X)$ is the numerator of the zeta function of $F.$\par
If $\ell =p,$ write $L_F(X)=\prod_{i}(1-\alpha_i X)$ and set $L_F^{nr}(X)=\prod_{v_p(\alpha_i)=0}(1-\alpha_i X).$ Then (see \cite{CRE} and
also \cite{GOS&SIN}):
$$\rm{Det}(1-\gamma X\mid_{T_{p}})=L_F^{nr} (X).$$\par
Now assume that $\ell$ does not divide the cardinal of $\Delta,$ then the above results are also valid character by character. More
precisely, if $\ell \not = p,$ we have:
$$\forall \chi \in \widehat{\Delta},\, \rm{Det}(1-\gamma X\mid_{T_{\ell}(\chi)})=L (X, \overline{\chi }).$$
If $\ell =p,$ for $\chi \in \widehat {\Delta},$ write $L(X,\chi)=\prod_i(1-\alpha_i(\chi )X)$ and set
$L^{nr}(X,\chi)=\prod_{v_p(\alpha_i(\chi )=0}(1-\alpha_i (\chi )X).$ Then:
$$\forall \chi \in \widehat{\Delta},\, \rm{Det}(1-\gamma X\mid_{T_{p}(\chi)})=L^{nr} (X, \overline{\chi }).$$\par
Now, let $\chi \in \widehat {\Delta},$ write:
$$L(X,\chi )=\prod_i(1-\alpha_i(\chi )X),$$
and set:
$$g(X,\chi)=\prod_{v_{\ell}(\alpha_i(\chi) -1)>0}(1-\alpha_i(\chi )X).$$
Set:
$$g(X)=\prod_{\chi \in \widehat{\Delta }}g(X, \chi ).$$
We also set:
$$\forall \chi \in \widehat {\Delta},\, H(X, \chi) =(1+X)^{\rm{deg}_X\, g(X,\chi )}g((1+X)^{-1},\chi ),$$
and:
$$H(X)=\prod_{\chi \in \widehat{\Delta }}H(X, \chi).$$\par
For $n\geq 0,$ set $F_n=\mathbb F_{q^{\ell ^n}}F,$ and let $A_n$ be the $\ell$-Sylow subgroup of $Cl^0(F_n).$ Let $F_{\infty}=\cup_{n\geq
0}F_n$ and let $A_{\infty}$ be the inductive limit of the $A_n,$ $n\geq 0.$ We set:
$$Y=\rm{Hom}(\mathbb Q_{\ell}/\mathbb Z_{\ell} , A_{\infty}).$$
Set $\Gamma=\rm{Gal}(F_{\infty}/F),$then $\gamma $ is a topological generator of $\Gamma \simeq \mathbb Z_{\ell}.$\par
\newtheorem{Lemma7}{Lemma}[section]
\begin{Lemma7} \label{Lemma7}
${}$\par
\noindent (1) For all $n\geq 0,$ we have an isomorphism of $\Delta$-modules:
$$\frac{Y}{(\gamma^{\ell^n}-1)Y}\simeq A_n.$$
(2) Assume $\mid \Delta \mid \not \equiv 0\pmod{\ell }.$Then, $\forall \chi \in \widehat{\Delta},$ $\forall n\geq 0,$ we have:
$$\frac{Y(\chi)}{(\gamma^{\ell^n}-1)Y(\chi)}\simeq A_n(\chi ).$$
\end{Lemma7}
\noindent{\sl Proof} We prove assertion (1), and note that (2) is a consequence of (1). Recall that $A_{\infty}$ is a divisible group (see
\cite{ROS}, Proposition 11.16). We start with the following exact sequence:
$$0\rightarrow A_n\rightarrow A_{\infty}\rightarrow A_{\infty}\rightarrow 0,$$
where the middle map is  the multiplication by $\gamma^{\ell^n}-1.$ We apply $\rm{Hom}(\mathbb Q_{\ell}/\mathbb Z_{\ell}, .)$ to this
sequence, we get:
$$0\rightarrow Y\rightarrow Y\rightarrow \rm{Ext}^1(\mathbb Q_{\ell}/\mathbb Z_{\ell}, A_n)\rightarrow 0.$$we also have the following exact
sequence:
$$0\rightarrow \mathbb Z_{\ell }\rightarrow \mathbb Q_{\ell}\rightarrow \frac{\mathbb Q_{\ell}}{\mathbb Z_{\ell}}\rightarrow 0.$$
We apply $\rm{Hom}(., A_n)$ to this last sequence, using the fact that:
$$\rm{Ext}^1(\mathbb Q_{\ell}, A_n)=\{ 0 \},$$
we get:
$$\rm{Hom}(\mathbb Z_{\ell },A_n)\simeq \rm{Ext}^1(\mathbb Q_{\ell}/\mathbb Z_{\ell} , A_n).$$
The Lemma follows. $\diamondsuit$\par
\newtheorem{Proposition4}[Lemma7]{Proposition}
\begin{Proposition4} \label{Proposition4}
${}$\par
 \noindent (1) Let $\Lambda =\mathbb Z_{\ell}[[X]]$ be the Iwasawa algebra of $\Gamma$ over $\mathbb Z_{\ell}$ where $X$ acts like $\gamma
-1.$ 
 Then $Y$ is a finitely generatyed $\Lambda$-module and a torsion $\Lambda$-module. The characteristic polynomial of the $\Lambda$-module $Y$
is equal to $H(X).$\par
\noindent (2) Assume that $\ell$ does not divide the cardinal of $\Delta.$ Let $\Lambda=W[[X]]$ be the Iwasawa algebra of $\Gamma$ over
$W=\mathbb Z_{\ell}[\mu_{\mid \Delta \mid}]$ where $X$ acts like $\gamma -1.$ Then, for $\chi \in \widehat{\Delta},$ $Y(\chi )$ is a finitely
generated $\Lambda$-module and a torsion $\Lambda$-module. The characteristic polynomial of tha $\Lambda$-module $Y$ is equal to
$H(X,\overline{\chi }).$
\end{Proposition4}
\noindent{\sl Proof} We prove (1), the proof of (2) is essentially similar. For all $n\geq 0,$ we set $\omega_n(X)=(1+X)^{\ell^n}-1.$ By
Lemma \ref{Lemma7}, we have:
$$\forall n\geq 0,\, \frac{Y}{\omega_n Y}\simeq A_n.$$
Therefore $Y$ is a finitely generated $\Lambda$-module and a torsion $\Lambda$-module. Let $r\in \mathbb N$ such that we have an isomorphism
of groups:
$$Y\simeq \mathbb Z_{\ell}^r.$$
Then, there exists a constant $\nu \in \mathbb Z,$ such that, for all $n$ sufficiently large:
$$\mid \frac{Y}{\omega_n Y} \mid =\ell^{rn+\nu}.$$
But, for all $n\geq 0,$ we have:
$$\mid A_n\mid =\ell^{v_{\ell}(L_{F_n}(1))}.$$
Therefore, there exists a constant $\nu'\in \mathbb Z$ such that, for all $n$ sufficiently large:
$$\mid A_n\mid = \ell ^{\rm{deg}_X H(X) n+\nu'}.$$
Thus: $r=\rm{deg}_X H(X).$ But let $V(X)$ be the characteristic polynomial of the $\Lambda$-module $Y.$ We know that $r=\rm{deg}_X V(X),$ and
we also know that $V(X)$ divides $(1+X)^{\rm{deg} L_F(X)}L_F((1+X)^{-1}).$ But $V(X)$ is a distinguished polynomial, thus $V(X)$ divides
$H(X).$ The Proposition follows. $\diamondsuit$\par
\newtheorem{Proposition5}[Lemma7]{Proposition}
\begin{Proposition5} \label{Proposition5}
${}$\par
\noindent (1) If $A_0$ is a cyclic $\mathbb Z_{\ell}$-module then $g(X)$ has simple roots.\par
\noindent (2) Assume that $\mid \Delta\mid \not \equiv 0\pmod{\ell }.$ Let $\chi \in \widehat{\Delta}.$ If $A_0(\chi )$ is a cyclic
$W$-module then $g(X,\overline{\chi })$ has simple roots.
\end{Proposition5}
\noindent{\sl Proof} We prove (1). By Nakayama's Lemma, $Y$ is pseudo-isomorphic to $\Lambda/H(X)\Lambda .$ But, by a result of Tate
(\cite{TAT}), we know that the action of $\gamma$ on $Y$ is semi-simple. This implies that $H(X)$ has simple roots. $\diamondsuit$\par
${}$\par
Let's give an application of this last Proposition.
\newtheorem{Proposition6}[Lemma7]{Proposition}
\begin{Proposition6} \label{Proposition6}
We assume that $q\geq 5.$ Let $E/\mathbb F_q(T)$ be a real quadratic field, i.e. $[E:\mathbb F_q(T)]=2$ and $\infty$ splits completely in
$E.$  If $O_E$ is a principal ideal domain then $L_E(X)$ has simple roots.
\end{Proposition6}
\noindent{\sl Proof} Let $g$ be the genus of $E$ and write:
$$L_E(X)=\prod_{i=1}^{2g}(1-\alpha_iX).$$
Let $K=\mathbb Q(\alpha_1 ,\cdots ,\alpha_{2g}),$ then $K$ is a CM-field. Let $\alpha \in \{
\alpha_1 ,\cdots \alpha_{2g} \}.$ Then:
$$(1-\alpha )(1-\overline{\alpha })\geq q+1-2\sqrt{q}>1.$$
Therefore:
$$N_{K/\mathbb Q}(1-\alpha )>1.$$
Thus $1-\alpha$ is not a unit of $K.$ Let $\infty_1$ and $\infty_2$ be tha places of $E$ above $\infty.$ Then $R(E)$ is a quotient of
$\mathbb Z (\infty_1-\infty_2)$ and we have an exact sequence:
$$0\rightarrow R(E) \rightarrow Cl^0(E)\rightarrow Cl(O_E)\rightarrow 0\, .$$
Therefore, if $O_E$ is a principal ideal domain then $Cl^0(E)$ is a cyclic group. It remains to apply Proposition \ref{Proposition5}.
$\diamondsuit$ \par
It is conjectured that there exists infinitely many real quadratic function fields $E/\mathbb F_q(T)$ such that $O_E$ is a principal ideal
domain. In view of this conjecture, it will be interesting to prove that there exists infinitely many real quadratic function fields
$E/\mathbb F_q(T)$ such that $L_E(X)$ has simple roots.\par
\section{A Conjecture of Goss}
${}$\par
Set $D_0=1$ and for $i\geq 1,$ $D_i=(T^{q^i}-T)D_{i-1}^q.$ The Carlitz exponential is defined by:
$$Exp(X)=\sum_{i\geq 0}\frac{X^{q^i}}{D_i}\in k[[X]].$$
Let $n\in\mathbb N,$ write $n=a_0+a_1q+\cdots +a_r q^r,$ where $a_0,\cdots , a_r\in \{ 0,\cdots ,q-1 \}.$  We set:
$$\Gamma_n=\prod_{i=0}^r D_i ^{a_i}.$$
The $i$th Bernoulli-Carlitz number, $B(i)\in k,$ is defined by:
$$\frac{X}{Exp(X)}=\sum_{i\geq 0}\frac{B(i)}{\Gamma_i}X^i.$$
Let $P$ be a prime of $A$ of degree $d$ and let $i\in\{ 1, \cdots , q^d-2\} ,$ $i\equiv 0\pmod{q-1}.$ We have the following result
(\cite{OKA}):
$$Cl(O_{K_P})_p(\omega_P^i)\not = \{ 0\} \Rightarrow B(i)\equiv 0\pmod{P}.$$
We fix an embedding of $\overline{\mathbb Q}$ in $\overline{\mathbb Q_p}.$ Let $i\in \{ 1,\cdots , q^d-2\}.$ Write:
$$L(X,\omega_P^i)=\prod_j(1-\alpha_j(i)X),$$
and set:
$$g(X,\omega_P^i)=\prod_{v_p(\alpha_j(i)-1)>0}(1-\alpha_j(i)X).$$\par
Let $i\in \mathbb N.$ We say that $i$ is a $q$-magic number if there exist  $c\in \{ 0, \cdots , q-2\} $ and an integer $n\in \mathbb N$
such that $i=cq^n+q^n-1.$\par
\newtheorem{Proposition7}{Proposition}[section]
\begin{Proposition7} \label{Proposition7}
Let $P$ be a prime of $A$ of degree $d.$ Let $i$ be a $q$-magic number, $1\leq i\leq q^d-2,$ $i\equiv 0\pmod{q-1}.$ Then $g(X,\omega_P^i)$
has simple roots.
\end{Proposition7}
\noindent{\sl Proof} We have $i=q^n -1$ for some integer $n,$ $1\leq n \leq d-1.$ By a result of Carlitz (\cite{GOS}, Lemma 8.22.4):
$$B(q^d-1-i)=\frac{(-1)^{d-n}}{L_{d-n}^{q^n}},$$
where $L_0=1$ and for $j\geq 1,$ $L_j=(T^{q^j}-T)L_{j-1}.$ Therefore:
$$Cl(O_{K_P})_p(\omega^{-i})=\{ 0\} .$$
It remains to aplly Proposition \ref{Proposition5}. $\diamondsuit$\par
${}$\par
In \cite{GOS}, David Goss makes the following conjecture:\par
\noindent let $P$ be a prime of degree $d$ and let $i$ be a $q$-magic number, $1\leq i\leq q^d-2.$ Then $\rm{deg}_Xg(X,\omega_P^i)\leq 1.$
\par
It is natural to ask if there exist primes $P$ and $q$-magic numbers $i$, $ 1\leq i\leq q^{\rm{deg}P}-2,$ such that
$\rm{deg}_Xg(X,\omega_P^i)\geq 1.$ This is the case.\par
\newtheorem{Proposition8}[Proposition7]{Proposition}
\begin{Proposition8} \label{Proposition8}
Let $c\in \{ 0,\cdots ,q-2\}.$ There exist infinitely many primes $P$ such that:
$$\prod_{n=1}^{\rm{deg}P-1}\beta(cq^n+q^n-1)\equiv 0\pmod{P}.$$
\end{Proposition8}
\noindent{\sl Proof} We prove this Proposition for $c\not =0.$ The proof for $c=0$ is very similar. If we apply the results in \cite{SHE}, we
get:
$$\forall n\geq 0,\, \rm{deg}_T\, \beta(cq^n+q^n-1)=n(c+1)q^n-\frac{q^{n+1}-q}{q-1}.$$
Let $S$ be the set of primes $P$ in $A$ such that:
$$\prod_{i=1}^{\rm{deg}P-1}\beta(cq^n+q^n-1)\equiv 0\pmod{P}.$$
Let's assume that $S$ is a finite set. We set:
$$D=\prod_{P\in S}\rm{deg}P,$$
and $D=1$ if $S=\emptyset.$ Note that:
$$\forall P\in S,\, q^D\equiv 1\pmod{q^{\rm{deg}P}-1}.$$
Therefore, since $\beta(c)=1,$ we have:
$$\forall P\in S,\, \beta(cq^D+q^D-1)\equiv 1\pmod{P}.$$
But $\rm{deg}_T\beta(cq^D+q^D-1)\geq 1,$ thus we can select a prime $Q$ of $A$ such that $\beta (cq^D+q^D-1)\equiv 0\pmod{Q}.$ Note that
$Q\not \in S.$ Set $d=\rm{deg}Q.$ Since $d$ does note divide $D,$ there exists  an integer $r,$ $1\leq r\leq d-1,$ such that $D\equiv
r\pmod{d}.$ Therefore:
$$\beta(cq^D+q^D-1)\equiv \beta (cq^r+q^r-1)\equiv 0\pmod{Q}.$$
But this implies that $Q\in S,$ which is a contradiction. $\diamondsuit$\par
${}$\par
Let $P$ be a prime of $A$ of degree $d.$ Let $J$ be the jacobian of $K_P,$ i.e. $J$ is the inductive limit of the $Cl^0(\mathbb F_{q^n}K_P),$
$n\geq 1.$ Set $\mathbb F_{q^{p^{\infty}}}=\cup_{n\geq 0}\mathbb F_{q^{p^n}}\subset \overline{\mathbb F_q},$ where $\overline{\mathbb F_q}$
is the algebraic closure of $\mathbb F_q$ in $\overline{k}.$ We consider the $\Delta=\rm{Gal}(K_P/k)$-module:
$${\cal{A}}_P=\frac{J[p]^{\rm{Gal}(\overline{\mathbb F_q}/\mathbb F_{q^{p^{\infty}}})}}{Cl^0(K_P)[p]}.$$
As a consequence of the results in section 4, we get:
\newtheorem{Proposition9}[Proposition7]{Proposition}
\begin{Proposition9} \label{Proposition9}
Let $W=\mathbb Z_{p}[\mu_{q^d-1}]$ and let $\chi \in \widehat{\Delta}.$ We have:
$$\rm{dim}_{\frac{W}{pW}}{\cal A}_P(\chi )=\rm{deg}_X g(X,\overline{\chi })-\rm{dim}_{\frac{W}{pW}}Cl^0(K_P)_p(\chi ).$$
\end{Proposition9}\par
Note that in general, by Proposition \ref{Proposition3}, we do not have ${\cal A}_P=\{ 0\}.$ But  Goss conjecture implies the
following:\par
\noindent let $P$ be a prime of $A$ of degree $d$ and let $i$ be a $q$-magic number, $1\leq i\leq q^d-2,$ then ${\cal A}_P(\omega_P^{-i})=\{
0\}.$\par
It would be interesting to prove (or find a counter-example) to this weak form of Goss
conjecture.\par

\end{document}